\documentclass[a4paper,11pt,reqno]{amsart}
\usepackage{mathrsfs}
 \usepackage[all]{xy}
\usepackage{amsmath,amssymb,amscd,bbm,amsthm,mathrsfs}
\usepackage{color,xcolor}
\usepackage{graphicx}
\usepackage{manfnt}
\newtheorem{thm}{Theorem}[section]
\newtheorem{lem}{Lemma}[section]

\newtheorem{con}{Conjecture}[section]
\newtheorem{pro}{Proposition}[section]

\begin{document}
\numberwithin{equation}{section}

 \title[Mixed curvature]{On mixed curvature for Hermitian manifolds}
\author{Kai Tang}
\address{Kai Tang. School of Mathematical Sciences, Zhejiang Normal University, Jinhua, Zhejiang, 321004, China} \email{{kaitang001@zjnu.edu.cn}}
\keywords{Hermitian manifold; Holomorphic sectional curvature; Chern Ricci curvature}
\thanks{\text{Foundation item:} Supported by National Natural Science
Foundation of China (Grant No.12001490).}
\begin{abstract}
In this paper, we consider {\em mixed curvature} $\mathcal{C}_{\alpha,\beta}$ for Hermitian manifolds, which is  a convex combination of the first Chern Ricci curvature and holomorphic sectional curvature introduced by Chu-Lee-Tam \cite{CLT}. We prove that if a compact Hermitian surface with constant mixed curvature $c$, then the Hermitian metric must be K\"{a}hler unless $c=0$ and $2\alpha+\beta=0$, which extends a previous result
by Apostolov-Davidov-Mu\v{s}karov \cite{ADM}. For the higher-dimensional case, we also partially classify compact locally conformal K\"{a}hler manifolds with constant mixed curvature. Lastly, we prove that if $\beta\geq0, \alpha(n+1)+2\beta>0$, then a compact Hermitian manifold with semi-positive but not identically zero mixed curvature has Kodaira dimension $-\infty$.
\end{abstract}

 \maketitle

\tableofcontents
\section{Introduction}
\subsection{Background}
For a Hermitian manifold $(M^{n},g,J)$, the Chern curvature $R$ is the curvature induced by the Chern connection $\nabla$ which is defined to be the connection such that $\nabla g=\nabla J=0$ with no $(1,1)$ components on the torsion. The first Chern Ricci curvature $Ric$ is defined by
\begin{align}
Ric=-\sqrt{-1}\partial\overline{\partial}\log \det g \nonumber \,\,,
\end{align}
where it is a $(1,1)$-form representing the first Chern class $c_{1}(M)$. The holomorphic sectional curvature $H$ is defined by
\begin{align}
H(X)=R(X,\overline{X},X,\overline{X})/|X|^{4} \nonumber \,\,,
\end{align}
for a (1, 0)-tangent vector $X\in T^{1,0}M$ (\cite{Zheng}). Given the perplexing relationship between the holomorphic sectional curvature and
the first Ricci curvature,  one can attempt to interpolate between these curvatures
by considering the following curvature constraint: In \cite{CLT} (see also \cite{BT}), Chu-Lee-Tam
introduced, for $\alpha,\beta\in\mathbb{R}$ the {\em mixed curvature}
 \begin{align}
\mathcal{C}_{\alpha,\beta}(X)=\frac{\alpha}{|X|^{2}_{g}}Ric(X,\bar{X})+\beta H(X)\nonumber \,\,.
 \end{align}
We always assume that $\alpha,\beta$ are both not $0$, otherwise, it makes no sense. $\mathcal{C}_{1,0}$ is the Ricci curvature, $\mathcal{C}_{0,1}$ is the holomorphic sectional curvature. When metric $g$ is K\"{a}hler, mixed curvature encompasses many interesting curvature conditions. $\mathcal{C}_{1,1}$ is the notion $Ric^{+}(X,\overline{X})$ introduced by Ni \cite{N2}. $\mathcal{C}_{1,-1}$ is the orthogonal Ricci curvature $Ric^{\perp}(X,\overline{X})$  introduced by Ni-Zheng \cite{NZ1}. $\mathcal{C}_{k-1,n-k}$ is closely related to the $k$-Ricci curvature $Ric_{k}$ introduced by Ni \cite{N1}.

Recently, there are many important results on compact K\"{a}hler manifolds with $\mathcal{C}_{\alpha,\beta}\geq 0$ or $\mathcal{C}_{\alpha,\beta}\leq 0$. It was proved by Yang \cite{Yang2018} that a compact K\"{a}hler manifold with $\mathcal{C}_{0,1}>0$ must be projective and rationally connected, confirming a conjecture of Yau \cite{Yau}. In \cite{N2}, Ni showed that it is also true if $Ric_{k}>0$ for some $1\leq k\leq n$. In \cite{Mat1,Mat2}, Matsumura established the structure theorems for a projective manifold with $\mathcal{C}_{0,1}\geq0$. In \cite{CLT, BT},
it was shown that any compact K\"{a}hler manifold with $\mathcal{C}_{\alpha,\beta}>0$, which constants satisfy $\alpha>0 $ and $3\alpha+2\beta\geq0$, must be projective and simply connected.  Zhang-Zhang \cite{ZZ} generalized Yang's result to quasi-positive
 case, confirming a conjecture of Yang affirmatively.  The author \cite{Tang2} showed the projectivity of compact K\"{a}hler manifolds with quasi-positive $\mathcal{C}_{\alpha,\beta}$. Very recently,
 Chu-Lee-Zhu \cite{CLZ} established a structure theorem for compact K\"{a}hler manifolds with $\mathcal{C}_{\alpha,\beta}\geq0$. For the non-positive case, Wu-Yau \cite{WuYau} confirmed a conjecture of Yau that a projective K\"{a}hler manifold with $\mathcal{C}_{0,1}<0$ must have ample canonical line bundle. Tosatti-Yang \cite{TosattiYang} was able to drop the projectivity assumption in Wu-Yau theorem. Chu-Lee-Tam \cite{CLT} proved that a compact K\"{a}hler manifold with $\mathcal{C}_{\alpha,\beta}<0$ and $\alpha\geq0, \beta\geq 0$, have ample canonical bundle. For more related works, we refer readers to \cite{BT2,DT,HW,N3,NZ2,Tang12,WuYau1,Yang20201,Zhang}.

Previous results have primarily focused on the study of K\"{a}hler cases. In this paper, we mainly investigate Hermitian manifolds.

\subsection{Constant mixed curvature}
The holomorphic sectional curvature plays a fundamental role in Hermitian geometry. Complete
K\"{a}hler manifolds with constant holomorphic sectional curvature $H$ are called complex space forms. They are quotients of complex projective space $\mathbb{CP}^{n}$, the complex Euclidean space $\mathbb{C}^{n}$, and the complex hyperbolic space $\mathbb{CH}^{n}$, equipped with (scaling of) the standard metrics. A long term folklore conjecture in non-K\"{a}hler geometry is the following:
\begin{con}
	Let $(M^{n},g)$ be a compact Hermitian manifold with $n\geq2$. Assume that $H=c$ where $c$ is a constant. If $c\neq0$, then $g$ is K\"{a}hler and if $c=0$, then $R=0$.
\end{con}
  Note that compact Chern flat manifolds were classified by Boothby \cite{Boothby} in 1958 as the set of
all compact quotients of complex Lie groups (equipped with left-invariant metrics compatible
with the complex structure). The compactness assumption in the above conjecture is necessary, as there are
 counterexamples in the noncompact case.

In complex dimension $n=2$, the  conjecture holds. In the $c\leq0$
case, it was solved by Balas-Gauduchon \cite{BG} (see also \cite{Balas}) in 1985 for Chern connection and by Sato-Sekigawa \cite{SS} in 1990 for Riemannian connection. The general case for $n=2$ (for both Chern connection and
Levi-Civita connection) was solved by Apostolov-Davidov-Mu\v{s}karov \cite{ADM} in 1996, as a corollary of their beautiful
classification theorem for compact self-dual Hermitian surfaces.

\begin{thm}[Apostolov-Davidov-Mu\v{s}karov] Any compact Hermitian surface with pointwise constant holomorphic sectional curvature with respect to the Levi-Civita connection or the Chern connection must be K\"{a}hler.
\end{thm}

For $n\geq3$, the conjecture is still largely open. In \cite{Tang1}, the author confirmed the conjecture under the additional assumption that
$g$ is Chern K\"{a}hler-like (namely $R$ obeys all K\"{a}hler symmetries). In \cite{CCN}, Chen-Chen-Nie proved the conjecture under the additional assumption that $g$ is locally conformally K\"{a}hler and $c\leq0$.  They also pointed out
 the necessity of the compactness assumption in the conjecture by explicit examples.  In \cite{ZhouZheng}, Zhou-Zheng  proved that any compact Hermitian threefold with vanishing
 {\em real bisectional curvature} must be Chern flat. Real bisectional curvature is a curvature notion
 introduced by  Yang-Zheng in \cite{YZ}. It is equivalent to holomorphic
 sectional curvature $H$ in strength when the metric is K\"{a}hler, but is slightly stronger than $H$
 when the metric is non-K\"{a}hler.  In \cite{RZ}, Rao-Zheng showed that the
 conjecture holds if the metric is Bismut K\"{a}hler-like, meaning that the curvature of the Bismut
 connection obeys all K\"{a}hler symmetries. Subsequently, Zhao-Zheng \cite{ZhaoZheng} proved a very beautiful theorem,  which states that all compact Strominger K\"{a}hler-like manifolds are pluriclosed. In \cite{BT3}, Broder and the author proved that pluriclosed metrics with vanishing holomorphic curvature on compact K\"{a}hler manifolds are K\"{a}hler. They also showed that Hermitian metrics with vanishing real bisectional curvature on compact complex manifolds in the Fujiki class $\mathcal{C}$ are K\"{a}hler.

On the other hand, if Chern Ricci curvature is non-zero constant $\lambda$, namely $Ric=\lambda g$, then $g$ is K\"{a}hler.
For these reasons, we naturally consider constant mixed curvature. We assume that $\alpha$ and $\beta$ are arbitrary fixed real numbers, and $\beta\neq 0$.  Mimicking the above folklore conjecture, we
propose the following:
\begin{con}
	Let $(M^{n},g)$ be a compact Hermitian manifold with $n\geq2$. Assume that mixed curvature $\mathcal{C}_{\alpha,\beta}=c$ where $c$ is a constant. If $c\neq0$, then $g$ is K\"{a}hler.
\end{con}

Note that if $\mathcal{C}_{\alpha,\beta}=0$, we may not be able to conclude that $R=0$. For example, isosceles Hopf manifold satisfies $\mathcal{C}_{\alpha,\beta}=0$ with $n\alpha+\beta=0$, but it does not imply $R=0$ which we will see in the section 3.

Using similar techniques as in \cite{ADM}, we prove confirm Conjecture 1.2 in the surface case:
\begin{thm}\label{1.2} If $(M^{2},g)$ is a compact Hermitian surface with constant mixed curvature $\mathcal{C}_{\alpha,\beta}=c$, then $g$ must be K\"{a}hler, unless $c=0$ and $2\alpha+\beta=0$; If $c=0$ and $2\alpha+\beta=0$,  the either $M^{2}$ is K\"{a}hler surface, or $M^{2}$ is an isosceles Hopf surface.
\end{thm}

In fact, we can  prove a stronger result, namely, assuming mixed curvature $\mathcal{C}_{\alpha,\beta}$ is pointwise constant $f$ where $f$ is a real-valued function $M^{2}$.
\begin{thm}\label{1.3} If $(M^{2},g)$ is a compact Hermitian surface with pointwise constant mixed curvature $\mathcal{C}_{\alpha,\beta}=f$. If $2\alpha+\beta\neq0$, then $g$ is K\"{a}hler.
\end{thm}

For higher-dimensional cases ($n\geq3$), we can also restrict ourselves to locally conformally K\"ahler metrics or Chern K\"{a}hler-like metrics. The second main result of this paper is the following theorem, which generalizes Chen-Chen-Nie's result \cite{CCN} to mixed curvature case.
\begin{thm}\label{1.4}
 Let $(M^{n},g)$ be a compact locally conformal K\"{a}hler manifold with constant mixed curvature $\mathcal{C}_{\alpha,\beta}=c$.
\begin{itemize}\label{1.4}
\item[(1)] If $c=0$ and $n\alpha+\beta\neq0$, then $g$ is K\"{a}hler.
\item[(2)] If $c<0$ and $[(n+2)\alpha+4\beta](n\alpha+\beta)\beta>0$, then $g$ is K\"{a}hler.
\item[(3)] If $c>0$ and $[(n+2)\alpha+4\beta](n\alpha+\beta)\beta<0$, then $g$ is K\"{a}hler.
\end{itemize}
\end{thm}

If the Hermitian metric is Chern K\"{a}hler-like, we can get similar conclusions even without the compactness condition.
\begin{thm}\label{1.5}
Let $(M^{n},g)$ be a  K\"{a}hler-like manifold with non-zero constant mixed curvature $\mathcal{C}_{\alpha,\beta}=c$. If $(n+2)\alpha+4\beta\neq0$, then $g$ is K\"{a}hler.
\end{thm}

For other connections, such as Strominger connection and Gauduchon connection, there have also been some new progress recently. Chen-Zheng \cite{CZ} proved that
any compact Hermitian surface with (pointwise)
constant Strominger holomorphic sectional curvature must be either K\"{a}hler or an isosceles Hopf surface with an admissible metric. Chen-Nie \cite{CN} showed that it is also true for any Hermitian surface with pointwise constant holomorphic sectional curvature with respect to a Gauduchon connection.
\subsection{Semi-positive mixed curvature}
The relationships between holomorphic sectional curvature and Ricci curvature, and the algebraic positivity of the (anti-)canonical line bundles, and some birational invariants of the ambient manifolds are still mysterious. In \cite{Yang2016}, Yang  proved that a compact Hermitian manifold with semi-positive
 but not identically zero holomorphic sectional curvature has Kodaira dimension $-\infty$. We extend Yang's result to mixed curvature.
\begin{thm}\label{1.6}
Let $(M^{n},g)$ be a  Hermitian manifold with  semi-positive but not identically zero mixed curvature $\mathcal{C}_{\alpha,\beta}$. If $\beta\geq0, \alpha(n+1)+2\beta>0$, then the Kodaira dimension $\kappa(M)=-\infty$.
\end{thm}

This paper is organized as follows. In section 2, we provide some basic knowledge which will be used in our proofs. In section 3, we prove  Theorem \ref{1.2}, \ref{1.3}, \ref{1.4} and \ref{1.5}. In section 4, we prove  Theorem \ref{1.6}.

\section{Preminaries}
Let $(M,g)$ be a Hermitian manifold of dimension $n$ and denote by $\omega$ the K\"{a}hler form associated with $g$. If $\omega$ is closed, that is, if $d\omega=0$, we call $g$ a K\"{a}hler metric. We will directly use $g$ to denote the  K\"{a}hler form, which should not cause any ambiguity. Denote by $\nabla$ Chern connection.

Fix any $p\in M$, let $\{e_{1},\cdot\cdot\cdot,e_{n}\}$ be a frame of (1,0)-tangent vector of $M$ in a neighborhood of $p$, with $
\{\varphi_{1},\cdot\cdot\cdot,\varphi_{n}\}$ being its dual coframe of (1,0)-form. We will also write $e=\,^{t}(e_{1},\cdot\cdot\cdot,e_{n})$ and
$\varphi=\,^{t}(\varphi_{1},\cdot\cdot\cdot,\varphi_{n})$ as colume vectors. Let $\langle\,,\,\rangle$ be the inner product given by Hermitian metric $g$, and extend it bi-linearly over $\mathbb{C}$. Denote $|\,\,|^{2}=\langle\,,\,\rangle$. Let $\theta$, $\tau$ and $\Theta$ be respectively the connection matrix, torsion colum vector and curvature matrix under $e$ for the Chern connection, then the structure equation and Bianchi identities are
\begin{align}
&d\varphi=-\,^{t}\theta \wedge \varphi+\tau \nonumber \\
&d\theta=\theta\wedge\theta+\Theta \nonumber \\
&d\tau=-\,^{t}\theta\wedge\tau+\,^{t}\Theta\wedge\varphi \nonumber \\
&d\Theta=\theta\wedge\Theta-\Theta\wedge\theta \nonumber
\end{align}
It is well know that the entries of matrix $\Theta$ are all (1,1)-form, while the entries of the column vector $\tau$ are all (2, 0)-forms, under any frame $e$. Denote
\begin{align}
\Theta_{kl}=R_{i\overline{j}k\overline{l}}\varphi_{i}\wedge \overline{\varphi}_{j} \nonumber
\end{align}

If we choose local holomorphic chart $(z_{1},\cdot\cdot\cdot,z^{n})$, we write
\begin{align}
\omega=\sqrt{-1}\sum_{i=1,j=1}^{n}g_{i\overline{j}}dz^{i}\wedge d\overline{z}^{j}. \nonumber
\end{align}
Then the curvature tensor $R=\{R_{i\overline{j}k\overline{l}}\}$ of the Chern connection is given by
\begin{align}
R_{i\overline{j}k\overline{l}}=-\frac{\partial^{2}g_{k\overline{l}}}{\partial z^{i}\partial \overline{z}^{j}}+g^{p\overline{q}}\frac{\partial g_{k\overline{q}}}{\partial z^{i}}\frac{\partial g_{p\overline{l}}}{\partial \overline{z}^{j}}   . \nonumber
\end{align}
There are four distinct traces of the Chern curvature tensor and, thus, four distinct Chern Ricci curvatures.  The \textit{first} and \textit{second Chern Ricci curvatures} $\text{Ric}_{i \bar{j}}^{(1)}  : =R_{i\overline{j}}=  \text{g}^{k \bar{\ell}} \text{R}_{i \bar{j} k \bar{\ell}}$ and  $\text{Ric}_{k \bar{\ell}}^{(2)}  : =  \text{g}^{i \bar{j}} \text{R}_{i \bar{j} k \bar{\ell}}$ trace to the same \textit{Chern scalar curvature} $u$.  The remaining \textit{third} and \textit{fourth Chern Ricci curvatures} $\text{Ric}_{i \bar{\ell}}^{(3)}  : =  \text{g}^{k \bar{j}} \text{R}_{i \bar{j} k \bar{\ell}}$ and $\text{Ric}_{k \bar{j}}^{(4)} : =  \text{g}^{i \bar{\ell}} \text{R}_{i \bar{j} k \bar{\ell}}$ trace to what we refer to as the \textit{altered Chern scalar curvature} $v$, i.e.,   \begin{eqnarray*}
 u\ : = \ \text{g}^{i \bar{j}} \text{Ric}_{i \bar{j}}^{(1)} \ = \ \text{g}^{k \bar{\ell}} \text{Ric}_{k \bar{\ell}}^{(2)}, \qquad v\ : = \ \text{g}^{i \bar{\ell}} \text{Ric}_{i \bar{\ell}}^{(3)} \ = \ \text{g}^{k \bar{j}} \text{Ric}_{k \bar{j}}^{(4)}.
\end{eqnarray*}
The first Chern Ricci curvature is a $(1,1)$-form representing the first Chern class $c_{1}(M)$. For convenience, we also denote the first, second, third, and fourth Ricci curvatures as $\rho^{(1)}$, $\rho^{(2)}$, $\rho^{(3)}$ and $\rho^{(4)}$ respectively.
 The components of the Chern torsion T are given by $\text{T}_{ij}^k = \text{g}^{k \bar{\ell}} (\partial_i \text{g}_{j \bar{\ell}} - \partial_j \text{g}_{i \bar{\ell}})$.  We write $\eta = \sum_i \eta_i dz^i = \sum_{i,k} \text{T}_{ik}^k dz^i$ for the \textit{torsion $(1,0)$-form}. It is well know that
 \begin{align}
\int_{M}(u-v)d\nu=\int_{M}|\eta|^{2}d\nu \nonumber
\end{align}
 Note that if $g$ is K\"{a}hler metric, then $R$ obeys all K\"{a}hler symmetries, four Chern Ricci curvatures are same and T$=0$.

By the conformal formula for Hermitian metric, the Hermitian $\widetilde{g}=e^{2F}g$ satisfies
 \begin{align}\label{2.01}
\widetilde{R}_{k\overline{l}i\overline{j}}=e^{2F}(R_{k\overline{l}i\overline{j}}-2g_{i\overline{j}}F_{k\overline{l}})
\end{align}
where $\widetilde{R}$ is the curvature tensor of $\widetilde{g}$ and $R$ is the curvature tensor of $g$.

For Hermitian surfaces $(M^{2},g)$, we list a few special formulas. By \cite[p113,(4.2)]{KO}, we get
 \begin{align}\label{2.2}
c_{1}^{2}(M)=\frac{1}{8\pi^{2}}[u^{2}-\langle Ric^{(1)},Ric^{(1)}\rangle]g^{2}
\end{align}
From \cite{BG} (also see \cite{G}), we have
\begin{align}\label{2.3}
&\rho^{(1)}+\rho^{(2)}-2 Re(\rho^{(3)})=(u-v)g \\
&\int_{M}\langle \rho^{(1)},\rho^{(2)}\rangle g^{2}=\int_{M}\langle \rho^{(1)}, Re(\rho^{(3)})\rangle g^{2}=
\int_{M}\langle\rho^{(1)},\rho^{(1)}\rangle g^{2}-\int_{M}u(u-v)g^{2}
\end{align}

Next let us recall the formula of Weyl curvature tensor for a Hermitian surface $(M^{2},g)$. One has $W=W^{+}+W^{-}$, and $g$ is said to be {\em self-dual} if $W^{-}=0$. Let $\{e_{1},e_{2}\}$ be a local unitary frame of $M$. Then
\begin{align}
\{e_{1}\wedge \overline{e}_{2}, e_{2}\wedge \overline{e}_{1}, \frac{1}{\sqrt{2}}(e_{1}\wedge\overline{e}_{1}-e_{2}\wedge\overline{e}_{2})\} \nonumber
\end{align}
form a basis of the complexification of $\Lambda^{2}_{-}$, the space on which the Hodge star operator is minus identity. We  call this basis the standard basis associated with the unitary frame $e$. In \cite{ADM}, it was proved that
\vspace{0.2cm}
\begin{lem}[\cite{ADM}]
Let $(M^{2},g)$ be a Hermitian surface, and $e$ a local unitary frame. Then under the standard basis associated with $e$, the components of $W^{-}$, the anti-self dual part of the Weyl tensor, are given by
\begin{align}
&W_{1}^{-}=R_{1\overline{2}1\overline{2}    }                                      \nonumber \\
&W_{2}^{-}= \frac{1}{\sqrt{2}}(R_{1\overline{2}2\overline{2}}+R_{2\overline{2}1\overline{2}}-
R_{1\overline{2}1\overline{1}}-R_{1\overline{1}1\overline{2}})                                          \nonumber \\
&W_{3}^{-}=\frac{1}{6}(R_{1\overline{1}1\overline{1}}+R_{2\overline{2}2\overline{2}}-
R_{1\overline{1}2\overline{2}}-R_{2\overline{2}1\overline{1}}-R_{1\overline{2}2\overline{1}}-R_{2\overline{1}1\overline{2}})                                          \nonumber
\end{align}
\end{lem}

\section{Constant mixed curvature}
\subsection{Hermitian surface} We will prove Theorem \ref{1.2} in this section. As in \cite{ADM}, the first step is to show that any Hermitian
surface with  constant mixed curvature must be self-dual:
\begin{pro}\label{3.1}
Let $(M^{2},g)$ be a  Hermitian surface. Assume that mixed curvature $\mathcal{C}_{\alpha,\beta}$ is pointwise constant, namely $\mathcal{C}_{\alpha,\beta}=f$ where $f$ is a real-valued function on $M^{2}$. Then $W_{-}=0$, namely, $M$ is self-dual.
\end{pro}
\noindent{{\em Proof.}} For any fixed $p\in M$, let $X=e_{1}$ and $Y=e_{2}$ be orthonormal $(1,0)$-vector, set $\mathcal{C}_{\alpha,\beta}=f(p)=c$. It follows from the average trick that
\begin{align*}
\frac{1}{vol(\mathbb{S}^{3})}\int_{ Z \in T^{1,0}_{p}M, | Z | =1 }\alpha Ric(Z,\overline{Z})+\beta H(Z) d\theta(Z)=(3\alpha+\beta)u+\beta v=6c \nonumber
\end{align*}
Hence
\begin{align}\label{3.1}
(3\alpha+2\beta)(R_{1\overline{1}1\overline{1}}+R_{2\overline{2}2\overline{2}})+(3\alpha+\beta)(R_{1\overline{1}2\overline{2}}+
R_{2\overline{2}1\overline{1}})+\beta(R_{1\overline{2}2\overline{1}}+R_{2\overline{1}1\overline{2}})=6c
\end{align}
By $\mathcal{C}_{\alpha,\beta}(X)=\mathcal{C}_{\alpha,\beta}(Y)=c$, we have
\begin{align}\label{3.2}
&\alpha(R_{1\overline{1}1\overline{1}}+R_{1\overline{1}2\overline{2}})+\beta R_{1\overline{1}1\overline{1}}=c \\
&\alpha(R_{2\overline{2}1\overline{1}}+R_{2\overline{2}2\overline{2}})+\beta R_{2\overline{2}2\overline{2}}=c
\end{align}
From (\ref{3.1}), (\ref{3.2}), (3.3) and $\beta\neq0$, we get
\begin{align}\label{3.4}
R_{1\overline{1}1\overline{1}}+R_{2\overline{2}2\overline{2}}-(R_{1\overline{1}2\overline{2}}+R_{2\overline{2}1\overline{1}}+
R_{1\overline{2}2\overline{1}}+R_{2\overline{1}1\overline{2}})=0
\end{align}
So $W_{-}^{3}=0$. Let $a,b\in\mathbb{C}$ and $|a|^{2}+|b|^{2}=1$. Then $aX+bY$ is unit vector. Thus, by carrying out the calculation, we obtain
\begin{align}\label{3.5}
c&=\alpha Ric(aX+bY,\overline{a}\overline{X}+\overline{b}\overline{Y})+\beta R(aX+bY,\overline{a}\overline{X}+\overline{b}\overline{Y},
aX+bY,\overline{a}\overline{X}+\overline{b}\overline{Y}) \nonumber\\
&=\alpha\{|a|^{2}Ric(X,\overline{X})+a\overline{b}Ric(X,\overline{Y})+b\overline{a}Ric(Y,\overline{X})+|b|^{2}Ric(Y,\overline{Y})\} \nonumber\\
&\,\,\,\,\,\,+\beta\{|a|^{4}R(X,\overline{X},X,\overline{X}+|b|^{4}R(Y,\overline{Y},Y,\overline{Y})) \nonumber \\
&\,\,\,\,\,\,+|a|^{2}|b|^{2}[R(X,\overline{X},Y,\overline{Y})+R(Y,\overline{Y},X,\overline{X})+R(X,\overline{Y},Y,\overline{X})+
R(Y,\overline{X},X,\overline{Y})]\nonumber \\
&\,\,\,\,\,\,+a^{2}\overline{b}^{2}R(X,\overline{Y},X,\overline{Y})+\overline{a}^{2}b^{2}R(Y,\overline{X},Y,\overline{X}) \nonumber\\
&\,\,\,\,\,\,+a^{2}\overline{a}\overline{b}[R(X,\overline{X},X,\overline{Y})+R(X,\overline{Y},X,\overline{X})]
+\overline{a}^{2}ab[R(Y,\overline{X},X,\overline{X})+R(X,\overline{X},Y,\overline{X})]   \nonumber \\
&\,\,\,\,\,\,+\overline{a}b^{2}\overline{b}[R(Y,\overline{Y},Y,\overline{X})+R(Y,\overline{X},Y,\overline{Y})]
+a\overline{b}^{2}b[R(Y,\overline{Y},X,\overline{Y})+R(X,\overline{Y},Y,\overline{Y})]  \}
\end{align}
In fact, by combining (\ref{3.4}), $\mathcal{C}_{\alpha,\beta}(X)=\mathcal{C}_{\alpha,\beta}(Y)=c$ and $|a|^{2}+|b|^{2}=1$, (\ref{3.5}) can be simplified to the following
\begin{align}\label{3.6}
0=&\alpha[a\overline{b}Ric(X,\overline{Y})+b\overline{a}Ric(Y,\overline{X})]+\beta
[a^{2}\overline{b}^{2}R(X,\overline{Y},X,\overline{Y})+\overline{a}^{2}b^{2}R(Y,\overline{X},Y,\overline{X})] \nonumber \\
&+\beta[(a^{2}\overline{a}\overline{b})A+(\overline{a}^{2}ab)B+(\overline{a}b^{2}\overline{b})C+(a\overline{b}^{2}b)D]
\end{align}
where we assume that
\begin{align}\label{3.7}
&A=R(X,\overline{X},X,\overline{Y})+R(X,\overline{Y},X,\overline{X});\,\,\,\,
B=R(Y,\overline{X},X,\overline{X})+R(X,\overline{X},Y,\overline{X}) \nonumber \\
&C=R(Y,\overline{Y},Y,\overline{X})+R(Y,\overline{X},Y,\overline{Y});\,\,\,\,\,\,\,\,
D=R(Y,\overline{Y},X,\overline{Y})+R(X,\overline{Y},Y,\overline{Y})
\end{align}
Now we choose special $a$ and $b$. Setting $a=\frac{1}{\sqrt{2}}, b=\frac{1}{\sqrt{2}}$, and $a=\frac{-1}{\sqrt{2}}, b=\frac{1}{\sqrt{2}}$. Then we get
\begin{align}\label{3.8}
R(X,\overline{Y},X,\overline{Y})+R(Y,\overline{X},Y,\overline{X})=0
\end{align}
Setting $a=\frac{\sqrt{2}}{2}, b=\frac{1-i}{2}$, and $a=\frac{-\sqrt{2}}{2}, b=\frac{1-i}{2}$. We have
\begin{align}\label{3.9}
R(X,\overline{Y},X,\overline{Y})-R(Y,\overline{X},Y,\overline{X})=0
\end{align}
Hence $R_{1\overline{2}1\overline{2}}=R(X,\overline{Y},X,Y)=0$. So $W_{-}^{1}=0$. Now  (\ref{3.6}) is further simplified to
\begin{align}\label{3.10}
0=&\alpha[a\overline{b}Ric(X,\overline{Y})+b\overline{a}Ric(Y,\overline{X})]
+\beta[(a^{2}\overline{a}\overline{b})A+(\overline{a}^{2}ab)B+(\overline{a}b^{2}\overline{b})C+(a\overline{b}^{2}b)D]
\end{align}
Setting $a=\frac{\sqrt{2}}{2}, b=-\frac{\sqrt{2}}{2}i$, and $a=\frac{\sqrt{3}}{2}, b=\frac{1}{2}i$, yields
\begin{align}
A+C=B+D \nonumber
\end{align}
Setting $a=-\frac{\sqrt{2}}{2}, b=\frac{\sqrt{2}}{2}$, and $a=\frac{\sqrt{3}}{2}, b=\frac{1}{2}$, yields
\begin{align}
A+B=C+D \nonumber
\end{align}
From these, we get
\begin{align}
R_{1\overline{2}2\overline{2}}+R_{2\overline{2}1\overline{2}}-R_{1\overline{2}1\overline{1}}-R_{1\overline{1}1\overline{2}}=D-A=0 \nonumber
\end{align}
namely, $W_{-}^{2}=0$. We conclude that $W_{-}=0$.
\qed

\vspace{0.2cm}
Next, we consider the case under conformal changes. We shall first establish a lemma.
\begin{lem}
Let $(M^{n},g)$ be a compact Hermitian manifold. Assume that mixed curvature $\mathcal{C}_{\alpha,\beta}$ is pointwise constant, namely $\mathcal{C}_{\alpha,\beta}=f$ where $f$ is a real-valued function on $M$. Then we have
\begin{align}\label{3.11}
[\alpha(n+2)+\beta]Ric^{(1)}+\beta Ric^{(2)}+2\beta Re (Ric^{(3)})=[2(n+1)f-\alpha u]g
\end{align}
\end{lem}

\noindent{{\em Proof.}} For $p\in M$, Let $X=X^{i}\frac{\partial}{\partial x^{i}}\in T_{p}^{(1,0)}M$. Since $\mathcal{C}_{\alpha,\beta}(X)=f(p)=c$ at $p$, then
\begin{align}
\alpha Ric(X,\overline{X})|X|^{2}+\beta R(X,\overline{X},X,\overline{X})=c|X|^{4}
\end{align}
Specifically, we have
\begin{align}
&\alpha X^{i}\overline{X}^{j}X^{k}\overline{X}^{l}(R_{i\bar{j}}g_{k\bar{l}}+R_{k\bar{j}}g_{i\bar{l}}+
R_{i\bar{l}}g_{k\bar{j}}+R_{k\bar{l}}g_{i\bar{j}}) \nonumber \\
&+\beta X^{i}\overline{X}^{j}X^{k}\overline{X}^{l}(R_{i\bar{j}k\bar{l}}+R_{k\bar{j}i\bar{l}}+R_{i\bar{l}k\bar{j}}+R_{k\bar{l}i\bar{j}})
=2cX^{i}\overline{X}^{j}X^{k}\overline{X}^{l}(g_{i\bar{j}}g_{k\bar{l}}+g_{i\bar{l}}g_{k\bar{j}})
\end{align}
Thus the above equality holds for any $X=X^{i}\frac{\partial}{\partial x^{i}}$ if and only if
\begin{align}\label{3.14}
\alpha (R_{i\bar{j}}g_{k\bar{l}}+R_{k\bar{j}}g_{i\bar{l}}+
R_{i\bar{l}}g_{k\bar{j}}+R_{k\bar{l}}g_{i\bar{j}})
+\beta(R_{i\bar{j}k\bar{l}}+R_{k\bar{j}i\bar{l}}+R_{i\bar{l}k\bar{j}}+R_{k\bar{l}i\bar{j}})
=2c(g_{i\bar{j}}g_{k\bar{l}}+g_{i\bar{l}}g_{k\bar{j}})
\end{align}
Taking the trace, we obtain (\ref{3.11}).                \qed

\begin{pro}\label{3.2}
Let $(M^{2},g)$ be a compact K\"{a}hler surface of constant holomorphic sectional curvature $H=c$. Let $\widetilde{g}=e^{2F}g$ is a conformal metric with
pointwise constant mixed curvature $\mathcal{C}_{\alpha,\beta}=f$, where $F$ is  a  smooth real-valued function on $M$.
\begin{itemize}
\item[(1)] If $2\alpha+\beta\neq0$, then $F$ is constant, namely, $\widetilde{g}$ is K\"{a}hler;
\item[(2)] If $f$ is constant, then $F$ is constant unless $f\equiv0$ and $2\alpha+\beta=0$.
\end{itemize}
\end{pro}
\noindent{{\em Proof.}} We denote the first, second  and third Ricci curvatures of $\widetilde{g}$ as $\widetilde{\rho}^{(1)},\widetilde{\rho}^{(2)},\widetilde{\rho}^{(3)}$  respectively. Let $u,v$ is scalar and altered  scalar curvature of $g$, and $\widetilde{u},\widetilde{v}$ is scalar and altered  scalar curvature of $\widetilde{g}$.

By (\ref{2.3}) and (\ref{3.11}), we have
\begin{align}
&\widetilde{\rho}^{(1)}+\widetilde{\rho}^{(2)}-2 Re (\widetilde{\rho}^{(3)})=(\widetilde{u}-\widetilde{v})\widetilde{g}  \\
&[4\alpha+\beta]\widetilde{\rho}^{(1)}+\beta\widetilde{\rho}^{(2)}+2\beta Re(\widetilde{\rho}^{(3)})=[6f-\alpha\widetilde{u}]\widetilde{g}
\end{align}
and  since $[3\alpha+\beta]\widetilde{u}+\beta\widetilde{v}=6f$, we obtain
\begin{align}\label{3.17}
2\alpha\widetilde{\rho}^{(1)}+2\beta Re (\widetilde{\rho}^{(3)})=(\alpha\widetilde{u}+\beta\widetilde{v})\widetilde{g}
\end{align}
By (2.4),
\begin{align}\label{3.18}
\int_{M}\langle\widetilde{\rho}^{(1)},Re (\widetilde{\rho}^{(3)})\rangle \widetilde{g}^{2}=\int_{M}\langle\widetilde{\rho}^{(1)},
\widetilde{\rho}^{(1)} \rangle\widetilde{g}^{2}-\int_{M}\widetilde{u}(\widetilde{u}-\widetilde{v})\widetilde{g}^{2}
\end{align}
Combining (\ref{3.17}) and (\ref{3.18}), we get
\begin{align}\label{3.18}
(\alpha+2\beta)\int_{M}\widetilde{u}^{2}\widetilde{g}^{2}=2(\alpha+\beta)\int_{M}\langle\widetilde{\rho}^{(1)},
\widetilde{\rho}^{(1)} \rangle\widetilde{g}^{2}+\beta\int_{M}(\widetilde{u}\widetilde{v})\widetilde{g}^{2}
\end{align}
Setting $c_{1}^{2}=\int_{M}c_{1}^{2}(M)$, according to (\ref{2.2}), we have
\begin{align}
c_{1}^{2}=\frac{1}{8\pi^{2}}[\int_{M}\widetilde{u}^{2}\widetilde{g}^{2}-\int_{M}\langle\widetilde{\rho}^{(1)},
\widetilde{\rho}^{(1)} \rangle\widetilde{g}^{2}]
\end{align}
Hence, we have
\begin{align}\label{3.21}
16\pi^{2}(\alpha+\beta)c_{1}^{2}=\alpha\int_{M}\widetilde{u}^{2}\widetilde{g}^{2}+\beta\int_{M}(\widetilde{u}\widetilde{v})\widetilde{g}^{2}
\end{align}
On the other hand, for K\"{a}hler metric $g$ with $H=c$, then $u=v=3c$. So we get
\begin{align}\label{3.22}
c_{1}^{2}=\frac{1}{16\pi^{2}}\int_{M}u^{2}g^{2}
\end{align}
From (\ref{3.21}) and (\ref{3.22}), we have
\begin{align}
(\alpha+\beta)\int_{M}u^{2}g^{2}=\alpha\int_{M}\widetilde{u}^{2}\widetilde{g}^{2}+\beta\int_{M}(\widetilde{u}\widetilde{v})\widetilde{g}^{2}
\end{align}
Since $\widetilde{g}=e^{2F}g$, then by (\ref{2.01}), we obtain
\begin{align}\label{3.24}
e^{2F}\widetilde{u}=u-4\Delta F;\,\,\,\,\,\,e^{2F}\widetilde{v}=v-2\Delta F
\end{align}
From these, we get
\begin{align}
(\alpha+\beta)\int_{M}u^{2}g^{2}=\alpha\int_{M}(u-4\Delta F)^{2}g^{2}+\beta \int_{M}(u-4\Delta F)(u-2\Delta F)g^{2}
\end{align}
Then
\begin{align}
(2\alpha+\beta)\int_{M}(\Delta F)^{2}g^{2}=0
\end{align}
where we use $\int_{M}\Delta F g^{2}=0$ and $u=3c$. This means that if $2\alpha+\beta\neq0$, then $\Delta F=0$. So $F$ is constant.

If $f$ is constant , we assume that $f=\widetilde{c}$. We only consider the case where $2\alpha+\beta=0$. We have
\begin{align}\label{3.27}
\alpha(\widetilde{u}-2\widetilde{v})=6\widetilde{c}
\end{align}
From (\ref{3.24}) and (\ref{3.27}), we get
\begin{align}
2\widetilde{c}e^{2F}=-\alpha c
\end{align}
If $\widetilde{c}\neq0$, then $F$ is constant.  \qed

\vspace{0.3cm}
Recall that Hopf manifold are compact complex manifolds whose universal cover is $\mathbb{C}^{n}\setminus\{0\}$. An isosceles Hopf manifold (see
\cite{CZ}) means $M^{n}_{\phi}=\mathbb{C}^{n}\setminus\{0\}/\langle \phi\rangle$ where
\begin{align}
\phi: (z_{1},z_{2},\cdot\cdot\cdot,z_{n})\rightarrow (a_{1}z_{1},a_{2}z_{2},\cdot\cdot\cdot,a_{n}z_{n})
\end{align}
with $0<|a_{1}|=|a_{2}|=\cdot\cdot\cdot=|a_{n}|<1$. The standard Hopf metric $g$ with K\"{a}hler form $\omega_{g}=\sqrt{-1}\frac{\partial\overline{\partial}|z|^{2}}{|z|^{2}}$ descends down to $M^{n}_{\phi}$. Write $h=|z|$, $e_{i}=h\frac{\partial}{\partial z_{i}}$, $\varphi_{i}=\frac{1}{h}dz_{i}$. Then $e$ becomes a unitary frame with $\varphi$ the
dual coframe. By the structure equations, we have
\begin{align}
\theta_{ij}=\frac{1}{h}(\bar{\partial}-\partial)h\delta_{ij}
\end{align}
Hence
\begin{align}
\Theta_{ij}=\frac{-2}{h^{2}}[\partial h\wedge(\overline{\partial}h)]\delta_{ij}+\frac{2}{h}(\partial\overline{\partial}h)\delta_{ij}
\end{align}
We get
\begin{align}
R_{i\overline{j}k\overline{l}}=-\frac{\overline{z}_{i}z_{j}}{|z|^{2}}\delta_{kl}+\delta_{ij}\delta_{kl}
\end{align}
where $h_{i}=\frac{\partial h}{\partial z^{i}}$. From this, the first Ricci curvature is
\begin{align}\label{3.33}
R_{i\overline{j}}=-n\frac{\overline{z}_{i}z_{j}}{|z|^{2}}+n\delta_{ij}
\end{align}
So we get scalar curvature $u=n^{2}-n$ and altered  scalar curvature $v=n-1$.
It implies
\begin{align}
(R_{i\bar{j}}g_{k\bar{l}}+R_{k\bar{j}}g_{i\bar{l}}+
R_{i\bar{l}}g_{k\bar{j}}+R_{k\bar{l}}g_{i\bar{j}})
-n(R_{i\bar{j}k\bar{l}}+R_{k\bar{j}i\bar{l}}+R_{i\bar{l}k\bar{j}}+R_{k\bar{l}i\bar{j}})
=0
\end{align}
So when $n\alpha+\beta=0$, the Hopf manifold have $\mathcal{C}_{\alpha,\beta}=0$. Note that the Hopf manifold is not Chern flat. We can also obtain the following lemma.

\begin{lem}\label{3.2.0}
Let $M^{2}$ be a Hopf surface with standard Hopf metric $g$. Let $\widetilde{g}=e^{2F}g$ is a conformal metric with
pointwise constant mixed curvature $\mathcal{C}_{\alpha,\beta}=f$, where $F$ is  a  smooth real-valued function on $M$.
Then  metric $\widetilde{g}$ does not exist, unless $f\equiv0$ and $2\alpha+\beta=0$.
\end{lem}
\noindent{{\em Proof.}} According to (\ref{3.21}) and (\ref{3.33}), it is easy to deduce that $c_{1}^{2} = 0$. From (\ref{3.21}), (\ref{3.24}) and $u=2v$, we have $\widetilde{u}=2\widetilde{v}$ and
\begin{align}
0=4\alpha\int_{M}(v-2\Delta F)^{2}g^{2}+2\beta\int_{M}(v-2\Delta F)^{2}g^{2}=2(2\alpha+\beta)\int_{M}(v-2\Delta F)^{2}g^{2}
\end{align}
it means that if $2\alpha+\beta\neq0$, then $\Delta F=\frac{v}{2}=\frac{1}{2}$, namely, $F$ is constant.  This of course will contradict with $v\neq0$.

On the other hand, it is easy to deduce that
\begin{align}
(2\alpha+\beta)\widetilde{v}=2f
\end{align}
When $2\alpha+\beta=0$, $f$ must be zero. \qed

\vspace{0.2cm}
Let $\widetilde{g}=e^{2F}g$ is a conformal metric with
pointwise constant mixed curvature $\mathcal{C}_{\alpha,\beta}=f$ on $M^{n}$. We denote  curvatures of $g$ as $R_{i\overline{j}k\overline{l}}$. Generally, from (\ref{2.01}) and (\ref{3.14}), we can obtain
\begin{align}\label{3.37}
&\alpha (R_{i\bar{j}}g_{k\bar{l}}+R_{k\bar{j}}g_{i\bar{l}}+
R_{i\bar{l}}g_{k\bar{j}}+R_{k\bar{l}}g_{i\bar{j}})
+\beta(R_{i\bar{j}k\bar{l}}+R_{k\bar{j}i\bar{l}}+R_{i\bar{l}k\bar{j}}+R_{k\bar{l}i\bar{j}}) \nonumber \\
&-2(n\alpha+\beta)[g_{i\overline{j}}F_{k\overline{l}}+g_{k\overline{j}}F_{i\overline{l}}+g_{i\overline{l}}F_{k\overline{j}}
+g_{k\overline{l}}F_{i\overline{j}}]
=2fe^{2F}(g_{i\bar{j}}g_{k\bar{l}}+g_{i\bar{l}}g_{k\bar{j}})
\end{align}

\vspace{0.2cm}
Now we are ready to prove the main result of this article, Theorem \ref{1.2}.

\vspace{0.2cm}
\noindent{\bf{{\em Proof of Theorem \ref{1.2}.}}}
Let $(M^{2},g)$ is a compact Hermitian surface with constant mixed curvature $\mathcal{C}_{\alpha,\beta}=c$. By Proposition \ref{3.1}, we know that
$M^{2}$ is is self-dual. Hence, by Theorem 1' in \cite{ADM}, $g$ is conformally to a metric $h$ on $M^{2}$, which is either (1) a complex space form, or (2) the non-flat, conformally flat K\"{a}hler metric, or (3) an isosceles Hopf surface. Write $g=e^{2F}h$.

For case (1), by Proposition \ref{3.2}, if $c\neq0$ or $2\alpha+\beta\neq0$, we get that $F$ is constant. Hence $g$ is K\"{a}hler.

For case (2), $M^{2}$ is a holomorphic bundle over a curve $C$ of genus at least 2, and $h$ is locally
a product metric where the factors have constant holomorphic sectional curvature -1 (for base curve) or 1 (for the fiber). Following \cite{ADM}, let $(z_{1},z_{2})$ be holomorphic coordinates, and K\"{a}hler form is
\begin{align}
\omega_{h}=2\sqrt{-1}\frac{dz^{1}\wedge d\overline{z}_{1}}{(1-|z_{1}|^{2})^{2}}+2\sqrt{-1}
\frac{dz^{2}\wedge d\overline{z}_{2}}{(1+|z_{2}|^{2})^{2}}=\omega_{1}+\omega_{2}
\end{align}
Let $e$ is unitary frame, such that $e_{i}$ is parallel to $\frac{\partial}{\partial z_{i}}$. Then under frame $e$, we have
\begin{align}
R_{1\overline{1}1\overline{1}}=-1;\,\,\,\,R_{2\overline{2}2\overline{2}}=1;
\end{align}
and others curvature is zero. So the Ricci curvature $R_{11}=-1$, $R_{2\overline{2}}=1$, others is zero. Denote $R_{i\overline{j}}=(-1)^{i}\delta_{ij}$. Here, let's assume that $g$ has pointwise constant mixed curvature $\mathcal{C}_{\alpha,\beta}=f$, where $f$ is  a real-valued function on $M^{2}$. According to formula (\ref{3.37}) , we can deduce that
\begin{align}\label{3.40}
-(\alpha+\beta)-2(2\alpha+\beta)F_{1\overline{1}}=fe^{2F};\,\,\,\, (\alpha+\beta)-2(2\alpha+\beta)F_{2\overline{2}}=fe^{2F};\,\,\,\,(2\alpha+\beta)F_{1\overline{2}}=0
\end{align}
If $2\alpha+\beta\neq0$. Setting $F_{1\overline{1}}=x$, $F_{2\overline{2}}=y$. So $y-x=\frac{\alpha+\beta}{2\alpha+\beta}$, $x_{\overline{2}}=F_{1\overline{1}\overline{2}}=0$ and $y_{\overline{1}}=F_{2\overline{2}\overline{1}}=0$. Then we have
\begin{align}
\Delta x=x_{1\overline{1}}=y_{1\overline{1}}=0;\,\,\,\, \Delta y=0
\end{align}
It implies that $\Delta F$ is constant. Hence $F$ is constant and $g$ is K\"{a}hler. Note that at this case we can only have  $f\equiv0$ and $\alpha+\beta=0$. If $2\alpha+\beta=0$, From (\ref{3.40}),
$f=0$ and $\alpha+\beta$=0, where contradicts with $2\alpha+\beta=0$.

For case (3), by Lemma \ref{3.2.0}, if $c\neq0$ or $2\alpha+\beta\neq0$,  the Hermitian metric $g$ does not exist on isosceles Hopf surface. This completes the proof of Theorem \ref{1.2}. \qed

\vspace{0.2cm}
\noindent{\bf{{\em Proof of Theorem \ref{1.3}.}}}
From the proof of Theorem \ref{1.2},  combining with Proposition \ref{3.2} and Lemma 3.2, we can easily draw the conclusion. \qed
\subsection{High-dimensional Hermitian manifold}
We first consider the case of compact locally conformal K\"{a}hler manifolds with constant mixed curvature. We also follow the similar steps as in \cite{CCN}. We first present the following lemma.
\begin{lem}
Let $(M^{n},g)$ be a compact locally conformal K\"{a}hler manifold with pointwise constant mixed curvature $\mathcal{C}_{\alpha,\beta}=f$.
\begin{itemize}
\item[(1)] If $f\equiv0$ and $n\alpha+\beta\neq0$, then $g$ is globally K\"{a}hler.
\item[(2)] if $f\leq0$ and $[(n+2)\alpha+4\beta](n\alpha+\beta)\beta>0$, then $g$ is globally K\"{a}hler.
\item[(3)] if $f\geq0$ and $[(n+2)\alpha+4\beta](n\alpha+\beta)\beta<0$, then $g$ is globally K\"{a}hler.
\end{itemize}
\end{lem}
\noindent{{\em Proof.}} We denote the first, second  and third Ricci curvatures of $g$ as $\rho^{(1)},\rho^{(2)},\rho^{(3)}$  respectively. Let $u,v$ is scalar and altered  scalar curvature of $g$. Since $g$ is locally conformal K\"{a}hler, then $\partial\partial^{\ast}g=\overline{\partial}\overline{\partial}^{\ast}g$ and
\begin{align}
\rho^{(1)}=\rho^{(3)}+\partial\partial^{\ast}g
\end{align}
where from (3.12) in \cite{CCN}. So $\rho^{(3)}$ is a real closed (1,1)-form. According to Proposition  3.2 in \cite{CCN}, we also have
\begin{align}
\rho^{(1)}-\rho^{(2)}=\frac{1}{n-1}[(v-u)g+n\partial\partial^{\ast}g]
\end{align}
By Lemma 3.1, we get
\begin{align}
&[\alpha(n+2)+\beta]\rho^{(1)}+\beta \rho^{(2)}+2\beta Re (\rho^{(3)})=[2(n+1)f-\alpha u]g;  \\
& [\alpha(n+1)+\beta]u+\beta v=n(n+1)f
\end{align}
From these, we obtain
\begin{align}
[\alpha(n+2)+4\beta]\rho^{(3)}=&\frac{1}{n-1}\{[(3n-2)(n+1)f-2(n\alpha+\beta)u]g \nonumber \\
&+[\beta(2-n)-\alpha(n+2)(n-1)]\partial\partial^{\ast}g\}
\end{align}\label{3.46}
Hence, We differentiate (\ref{3.46}) to obtain
\begin{align}\label{3.47}
d\{[(3n-2)(n+1)f-2(n\alpha+\beta)u]g\}=0
\end{align}
Let $\eta=(3n-2)(n+1)f-2(n\alpha+\beta)u$. Then we will prove that $\eta$ is either everywhere nonzero or $\eta\equiv0$.

(1) If $\eta$ is everywhere nonzero, then $g$ is globally K\"{a}hler.

(2) If $\eta=0$ at some point $p$ on $M$, then $\eta\equiv0$ on $M$ (see \cite{CCN}). Therefore, if $\eta\equiv0$ we have the following case:

(a) If $f\equiv0$ and $n\alpha+\beta\neq0$, then $u=v=0$. It implies that $g$ is balanced. Since $g$ is also locally conformal K\"{a}hler, we get $g$ is K\"{a}hler.

(b) if $f\leq0$ and $[(n+2)\alpha+4\beta](n\alpha+\beta)\beta>0$. We only consider the case where there is at least one point $p$ such that $f(p)<0$. Since
\begin{align}
u=\frac{(3n-2)(n+1)f}{2(n\alpha+\beta)},\,\,\,\,\,\,\,\, v=\frac{-(n+1)[(n+2)(n-1)\alpha+(n-2)\beta]f}{2\beta(n\alpha+\beta)}
\end{align}
Hence, we get
\begin{align}
v-u=\frac{-(n+1)(n-1)[(n+2)\alpha+4\beta]f}{2(n\alpha+\beta)\beta}
\end{align}
So $\int_{M}(v-u)g^{n}>0$, which is contradiction with a well-known result that $\int_{M}(v-u)g^{n}\leq0$. If $f\geq0$ and $[(n+2)\alpha+4\beta](n\alpha+\beta)\beta<0$, this would also lead to a contradiction.

In this way, we have completed the proof.     \qed

\vspace{0.2cm}
\noindent{\bf{{\em Proof of Theorem \ref{1.4}.}}}
Let $g$ be a compact locally conformal K\"{a}hler manifold with constant mixed curvature $\mathcal{C}_{\alpha,\beta}=c$. Under the assumptions of Lemma 3.3, we know that $g$ is globally K\"{a}hler.
We assume that $h$ is K\"{a}hler metric where $dh=0$, and $g=e^{2F}h$. Setting
\begin{align}\label{3.50}
A=[(3n-2)(n+1)c-2(n\alpha+\beta)u]e^{2F}
\end{align}
By (\ref{3.47}), we get  that $d(Ah)=0$. So $A$ is constant. Since $h$ is K\"{a}hler, we get
\begin{align}\label{3.51}
e^{2F}u=e^{2F}v+2(1-n)\Delta_{h} F
\end{align}
We also have that
\begin{align}\label{3.52}
(\alpha(n+1)+\beta)u+\beta v=n(n+1)c
\end{align}
From (\ref{3.50}), (\ref{3.51}) and (\ref{3.52}), we can easily deduce that
\begin{align}\label{3.53}
ce^{2F}(n+1)(n-1)[(n+2)\alpha+4\beta]=[\alpha(n+1)+2\beta]A+4\beta(1-n)(n\alpha+\beta)\Delta_{h} F
\end{align}

(1) If $c=0$ and $n\alpha+\beta\neq0$, By integrating (\ref{3.53}) with respect to the metric $h$, we obtain that $F$ is constant. It implies that $g$ is K\"{a}hler.

(2) If $c<0$, $[(n+2)\alpha+4\beta](n\alpha+\beta)\beta>0$, we will show that $F$ is constant. Setting
\begin{align}
B=[\alpha(n+1)+2\beta]A
\end{align}
By integrating (\ref{3.53}) with respect to the metric $h$, we have
\begin{align}\label{3.55}
c(n+1)(n-1)[(n+2)\alpha+4\beta]\int_{M}e^{2F}h^{n}=B\cdot Vol(M)
\end{align}
On the other hand, we also have that
\begin{align}\label{3.56}
\Delta_{h} F=-\frac{1}{2}e^{2F}(\Delta_{h} e^{-2F})+2e^{2F}|\partial F|^{2}_{g}
\end{align}
Substituting equation (\ref{3.56}) into equation (\ref{3.53}), and then integrating it, we obtain
\begin{align}\label{3.57}
B\int_{M}e^{-2F}h^{n}=&c(n+1)(n-1)[(n+2)\alpha+4\beta]Vol(M) \nonumber \\
&+8\beta(n-1)(n\alpha+\beta)\int_{M}|\partial F|^{2}_{g}h^{n}
\end{align}
By (\ref{3.55}) and (\ref{3.56}) and Cauchy-Schwarz inequality, we have
\begin{align}\label{3.58}
(Vol(M))^{2}\leq \int_{M}e^{2F}h^{n}\int_{M}e^{-2F}h^{n}\leq (Vol(M))^{2} \nonumber
\end{align}
The equality holds if and only if $F$ is constant. Therefore that $g$ is K\"{a}hler.

(3) If $c>0$ and $[(n+2)\alpha+4\beta](n\alpha+\beta)\beta<0$, the proof is similar to case (2). \qed

\vspace{0.2cm}
Now we will prove the K\"{a}hler-like case.

\noindent{\bf{{\em Proof of Theorem \ref{1.5}.}}} This proof is relatively simple. Let $(M^{n},g)$ be a  K\"{a}hler-like manifold with constant mixed curvature $\mathcal{C}_{\alpha,\beta}=c$. Then $Ric^{(1)}=Ric^{(2)}=Ric^{(3)}$, $u=v$ and by Lemma 3.1, we have
\begin{align}
&[(n+2)\alpha+4\beta]Ric^{(1)}=[2(n+1)c-\alpha u]g \nonumber \\
&[(n+1)\alpha+2\beta]u=n(n+1)c \nonumber
\end{align}
When $(n+2)\alpha+4\beta\neq0$, we get that
\begin{align}
[(n+1)\alpha+2\beta]Ric^{(1)}=(n+1)cg  \nonumber
\end{align}
If $c\neq0$, then $g$ is K\"{a}hler.       \qed

\section{Semi-positive mixed curvature}
In this section, we provide the proof of Theorem \ref{1.6}.

\vspace{0.2cm}
\noindent{\bf{{\em Proof of Theorem \ref{1.6}.}}}
Let $\text{g}_{\text{G}} = f_0^{\frac{1}{n-1}} \text{g}$ be the Gauduchon metric (i.e., g$_{\text{G}}$ satisfies $\partial \bar{\partial} \omega_{\text{G}}^{n-1}=0$) in the conformal class of g, where $f_0 \in \mathcal{C}^{\infty}(X,\mathbf{R})$ is a strictly positive function (such a metric always exists \cite{G1}). Let $\omega_{\text{G}}(\cdot, \cdot) : = \text{g}_{\text{G}}(\text{J} \cdot, \cdot)$. Let $S_{g_{G}}$ and $\widetilde{S}_{g_{G}}$ is scalar and altered  scalar curvature of $g_{G}$, $S_{g}$ and $\widetilde{S}_{g}$ is scalar and altered  scalar curvature of $g$ .

Then \begin{eqnarray*}
\int_X \text{S}_{\text{g}_{\text{G}}} \omega_{\text{G}}^n &=& n \int_X \text{Ric}_{\text{g}_{\text{G}}}^{(1)} \wedge \omega_{\text{G}}^{n-1} \\
&=& n \int_X \text{Ric}_{\text{g}}^{(1)} \wedge \omega_{\text{G}}^{n-1} \ = \  n \int_X f_0 \text{Ric}_{\text{g}}^{(1)} \wedge \omega_{\text{g}}^{n-1} \ = \  \int_X f_0 \text{S}_{\text{g}} \omega_{\text{g}}^n.
\end{eqnarray*}
A similar computation (see, e.g., \cite{Yang2016}) then gives $\int_X \widetilde{\text{S}}_{\text{g}_{\text{G}}} \omega_{\text{G}}^n = \int_X f_0 \widetilde{\text{S}}_{\text{g}} \omega_{\text{g}}^n$.  Let $\eta_{\text{G}}$ denote the torsion $(1,0)$-form of $\text{g}_{\text{G}}$. We then compute  \begin{eqnarray*}
(\alpha(n+1)+2\beta)\int_X \text{S}_{\text{g}_{\text{G}}} \omega_{\text{G}}^n &=& \int_X \left((\alpha(n+1)+\beta) \text{S}_{\text{g}_{\text{G}}} + \beta\widetilde{\text{S}}_{\text{g}_{\text{G}}} \right) \omega_{\text{G}}^n + \int_X \beta\left( \text{S}_{\text{g}_{\text{G}}} - \widetilde{\text{S}}_{\text{g}_{\text{G}}} \right) \omega_{\text{G}}^n  \\
&=& \int_X f_0 \left((\alpha(n+1)+\beta) \text{S}_{\text{g}} +\beta \widetilde{\text{S}}_{\text{g}} \right) \omega_g^n +\beta\int_X | \eta_{\text{G}} |^2_{\text{g}_{\text{G}}} \omega_{\text{G}}^n  \\
\end{eqnarray*}
It follows from the average trick that
\begin{align*}
\frac{1}{vol(\mathbb{S}^{2n-1})}\int_{ Z \in T^{1,0}_{p}M, | Z | =1 }\alpha Ric_{g}(Z,\overline{Z})+\beta H_{g}(Z) d\theta(Z)=\frac{((n+1)\alpha+\beta)S_{g}+\beta \widetilde{S}_{g}}{n(n+1)} \nonumber
\end{align*}
If $g$ has semi-positive but not identically zero mixed curvature $\mathcal{C}_{\alpha,\beta}$ and $\beta\geq0, \alpha(n+1)+2\beta>0$,
then we have
\begin{eqnarray*}
\int_X \text{S}_{\text{g}_{\text{G}}} \omega_{\text{G}}^n>0
\end{eqnarray*}
It means that the Gauduchon metric  $\text{g}_{\text{G}}$ has positive total scalar curvature. By Theorem 1.3 in \cite{Yang2019}, we know that $M$ admits a Hermitian metric with positive scalar curvature, which it implies that Kodaira dimension $\kappa(M)=-\infty$. \qed
\vspace{0.3cm}

\vspace{0.5cm}
\noindent\textbf{Acknowledgement.} The author is grateful to Professor Fangyang Zheng for constant encouragement and support.

\end{document}